\documentclass[preprint]{elsarticle}
\usepackage[T1]{fontenc}
\usepackage[latin9]{inputenc}
\usepackage[letterpaper]{geometry}
\usepackage{float}
\usepackage{amsthm}
\usepackage{relsize}
\usepackage{graphicx}
\usepackage{amssymb}
\usepackage{amsmath}
\usepackage[american,english]{babel}
\makeatletter

\theoremstyle{plain}
\newtheorem{thm}{Theorem}
\newtheorem{prop}[thm]{Proposition}
\newtheorem{theorem}{Theorem}

\newtheorem{remark}[theorem]{Remark}

\journal{Annales de l'Institut Henri Poincar\'{e} (B) Probabilit\'{e}s et Statistiques}

\begin{document}
\begin{frontmatter}

\title{On the density of the sum of two independent Student t-random vectors}

\author[ku]{C. Berg}
\ead{berg@math.ku.dk}

\author[mlv]{C. Vignat\corref{cor1}}
\ead{vignat@univ-mlv.fr}

\cortext[cor1]{Corresponding author}

\address[ku]{Department of Mathematics, University of Copenhagen, Universitetsparken 5, DK-2100, Copenhagen, Denmark}

\address[mlv]{Institut Gaspard Monge, Université de Marne la Vall\'{e}e, France}

\begin{abstract}
In this paper, we find an expression for the density of the sum of two independent $d-$dimensional Student $t-$random vectors $\mathbf{X}$ and $\mathbf{Y}$ with arbitrary degrees of freedom. As a byproduct we also obtain an expression for the density of the sum $\mathbf{N}+\mathbf{X}$, where $\mathbf{N}$ is normal and $\mathbf{X}$ is an independent Student $t-$vector.
In both cases the density is given as an infinite series
\[
\sum_{n=0}^{\infty} c_{n}f_{n}
\]
where $f_{n}$ is a sequence of probability densities on $\mathbb{R}^{d}$ and $\left( c_{n} \right)$ is a sequence of positive numbers of sum $1$, i.e. the distribution of a non-negative integer-valued random variable $C$, which turns out to be infinitely divisible for $d=1$ and $d=2.$ When $d=1$ and the degrees of freedom of the Student variables are equal, we recover an old result of Ruben.
\end{abstract}

\end{frontmatter}

\section{Introduction}
The Student $t-$distributions plays an important role in statistics for the analysis of normal samples. One of the major properties of these distributions is their heavy-tail behaviour; the same behaviour is shown by the famous family of stable distributions. 
The Student t-distributions show the advantage of having an explicit expression, what allows for example to write down explicitly their likelihood function \cite{Bening}; this is not the case for the stable distributions, except some well-known particular cases.

However, except in the Gaussian and Cauchy case, the Student t-distributions  show the drawback of  not being stable by convolution. 
Thus the study of their behavior  by convolution is of high interest. This study has been first addressed in \cite{Walker} and more recently in \cite{Witkovsky, BergVignat1}, but concerns only the case where the degrees of freedom is an odd integer.
The study of Student t-processes has recently given rise to a series of studies \cite{Cufaro, BergVignat2, Leonenko}.

The $d-$dimensional Student $t-$density with
$2\nu>0$ degrees of freedom is 
\[
f_{\nu}\left(\mathbf{t}\right)=A_{\nu,d}\left(1+\Vert \mathbf{t}\Vert ^{2}\right)^{-\nu-\frac{d}{2}},\,\, \mathbf{t}\in\mathbb{R}^{d};\,\, A_{\nu,d}=\frac{\Gamma\left(\nu+\frac{d}{2}\right)}{\Gamma\left(\nu\right)\pi^{\frac{d}{2}}}.
\]

In the following, we provide an expression for the density of the sum 
\[
\mathbf{Y}=\mathbf{T}_{1}+\mathbf{T}_{2}
\]
where $\mathbf{T}_{1}$ and $\mathbf{T}_{2}$ are $d-$variate independent $t-$distributed random
variables with respective degrees of freedom $2\nu$ and $2\mu$. When $\mu = \nu$ and $d=1,$ this equals an expression found by Ruben, cf. \cite[formula 3.4]{Ruben}.
As a byproduct, we also obtain an expression for the density of
\[
\mathbf{Z}=\mathbf{X}+\mathbf{T},
\]
where $\mathbf{X}$ is standard Gaussian.
Moreover, we extend this result to the case of $\tilde{\mathbf{Z}}=a_{1}\mathbf{X}+a_{2}\mathbf{T}$ with $a_{1},a_{2}>0$. In all cases, our expression is an infinite convex combination of a sequence of $d-$variate densities, which we describe next.

Let us introduce the family of densities
\begin{equation}
g_{k,\sigma}(\mathbf{z})
=\frac{1}{\sigma^{d}}g_{k,1}\left(\frac{\mathbf{z}}{\sigma}\right)
=\frac{\Gamma(\frac{d}{2})}{\Gamma\left(k+\frac{d}{2}\right)\left(\sigma\sqrt{2\pi}\right)^{d}}\left(\frac{\Vert \mathbf{z} \Vert ^{2}}{2\sigma^{2}}\right)^{k}\exp\left(-\frac{\Vert \mathbf{z}\Vert^{2}}{2\sigma^{2}}\right),\,\, \mathbf{z}\in\mathbb{R}^{d}\label{eq:gk}
\end{equation}
with $k\in\mathbb{N}.$ We remark that $g_{0,\sigma}$ is the $d-$variate  Gaussian
density with zero mean and covariance matrix $\sigma^{2}I_{d};$ for $k\ne0$, by
\cite[Th. 2.9]{Fang}, $g_{k,\sigma}$ is the density of $\sigma\Vert\mathbf{X}_{2k+d}\Vert \mathbf{U}_{d}$ where $\mathbf{X}_{2k+d}\in\mathbb{R}^{2k+d}$
is a Gaussian vector with identity covariance matrix and $\mathbf{U}_{d}$ is uniformly distributed on the unit sphere $\mathcal{S}_d$ of $\mathbb{R}^d$ and independent of $\mathbf{X}_{2k+d}.$

Let us also introduce the family of densities
\begin{equation}
\phi_{n,\eta}\left(\mathbf{x}\right)=\frac{\Gamma\left( \frac{d}{2} \right)}{\pi^{\frac{d}{2}} B\left(\eta,n+\frac{d}{2}\right)}\frac{\Vert \mathbf{x}\Vert ^{2n}}{\left(1+\Vert \mathbf{x}\Vert^{2}\right)^{\eta+n+\frac{d}{2}}},\,\,\mathbf{x}\in\mathbb{R}^{d},\label{eq:sqrtF}
\end{equation}
where $\eta>0.$
We note that, with the same notation as above, $\phi_{n,\eta}$ is
the density of $\Vert\mathbf{T}_{2n+d}\Vert \mathbf{U}_{d},$ where $\mathbf{T}_{2n+d}\in\mathbb{R}^{2n+d}$
follows a multivariate Student t-distribution with $2\eta$ degrees
of freedom
\[
\frac{\Gamma\left(\eta+n+\frac{d}{2}\right)}{\Gamma\left(\eta\right)\pi^{n+\frac{d}{2}}}\left(1+\Vert \mathbf{x}\Vert^{2} \right)^{-\eta-n-\frac{d}{2}},\,\,\mathbf{x}\in \mathbb{R}^{2n+d}.
\]

At last, we recall that a discrete random variable $B$ follows a
negative binomial distribution with parameters $r$ and $p$ (denoted
as $NB\left(r,p\right)$) if
\[
\Pr\left\{ B=k\right\} =
{k+r-1 \choose k}p^{r}\left(1-p\right)^{k},\,\,r>0,\,\,0<p<1,\,\,k=0,1\dots
\]
$NB$ distributions with non-integer parameter $r$ are also known
as Pólya distributions.

\section{Results}

The distribution of the sum of two independent Student $t-$vectors
\[
\mathbf{Y}=\mathbf{T}_{1}+\mathbf{T}_{2}
\]
can be deduced by subordination from the distribution of $\mathbf{Z}=\mathbf{X}+\mathbf{T}$ where $\mathbf{X}$ is Gaussian 
: if $\tilde{\gamma}_{\nu}$ denotes the inverse Gamma
measure with scale parameter $\frac{1}{4}$ and shape parameter $\nu>0$
\[
d\tilde{\gamma}_{\nu}\left(t\right)=\frac{1}{2^{2\nu}\Gamma\left(\nu\right)}\exp\left(-\frac{1}{4t}\right)t^{-\nu-1}dt,\,\,t>0,
\]
and if $g_{t}\left(\mathbf{x}\right)$ denotes the $d-$variate Gaussian semi-group of
normal densities
\[
g_{t}\left(\mathbf{x}\right)=\left( 4\pi t \right)^{-\frac{d}{2}} \exp\left(-\frac{\Vert \mathbf{x}\Vert^{2}}{4t}\right),\,\,t>0,\,\,\mathbf{x}\in \mathbb{R}^{d},
\]
then the Student $t-$density can be obtained as the mixture 
\[
f_{\nu}\left(\mathbf{x}\right)=\int_{0}^{+\infty}g_{t}\left(\mathbf{x}\right)d\tilde{\gamma}_{\nu}\left(t\right).\]
Therefore, we begin by finding an expression for the density of $\mathbf{Z}$.

\subsection{The sum of a Gaussian and a Student $t-$vector}

Our first result states as follows.
\begin{thm}
\label{thm:Student+Gaussian}
The density of $\mathbf{Z}=\mathbf{X}+\mathbf{T}$ can be written 
\begin{equation}
f_{\mathbf{Z}}(\mathbf{z})=\sum_{k=0}^{+\infty}\alpha_{k}^{(\nu,\gamma)}g_{k,\sigma}(\mathbf{z})
\label{eq:fzmultivariate}
\end{equation}
where the densities $g_{k,\sigma}$ are defined in (\ref{eq:gk}), $\alpha_{k}^{(\nu,\gamma)}$
are positive coefficients that sum to $1$ and given by 
\begin{equation}
\alpha_{k}^{(\nu,\gamma)}=\frac{\Gamma(k+\frac{d}{2})}{k!\Gamma\left(\frac{d}{2}\right)\Gamma(\nu)}\gamma^{2k}\int_{0}^{+\infty}\exp(-a)a^{\nu+\frac{d}{2}-1}\left(a+\gamma^{2}\right)^{-k-\frac{d}{2}}da,\,\,\, k\ge0
\label{eq:alphakmultivariate}
\end{equation}
and $\gamma$ is defined as 
\begin{equation}
\gamma=\frac{1}{\sigma\sqrt{2}}.
\label{eq:gamma}
\end{equation}
\end{thm}

\begin{prop}
\label{pro:extensionStudentGaussian}
The same result as in Theorem
\ref{thm:Student+Gaussian} holds for $\tilde{\mathbf{Z}}=a_{1}\mathbf{X}+a_{2}\mathbf{T},\,\,\,a_{1}\in\mathbb{R}^{+},\,a_{2}\in\mathbb{R}^{+}$
replacing $\sigma$ by $ a_{1}\sigma$ and $\gamma$
by $\frac{a_{2}}{\sigma\sqrt{2} a_{1}}.$
\end{prop}
Note that $a_1\mathbf{X}$ and $a_2\mathbf{T}$ have the same distributions as $\left(-a_1\mathbf{X}\right)$ and $\left(-a_2\mathbf{T}\right)$ respectively, so that we only need to consider the case $a_1>0,\,a_2>0$ in the above Proposition.

\subsection{The sum of two Student $t-$vectors}

Starting from the distribution (\ref{eq:fzmultivariate}) we obtain the
following
\begin{thm}
\label{thm:Student+Student}
The distribution of $\mathbf{Y}=\mathbf{T}_{1}+\mathbf{T}_{2},$ where
$\mathbf{T}_{1}$ and $\mathbf{T}_{2}$ are independent $t-$distributed with respective
degrees of freedom $2\nu>0$ and $2\mu>0,$ can be written as the infinite convex combination
\begin{equation}
\label{fT1T2}
f_{\mathbf{Y}}\left(\mathbf{x}\right)=\sum_{n=0}^{+\infty}c_{n}^{\left(\nu,\mu\right)}\phi_{n,\nu+\mu}\left(\mathbf{x}\right),
\end{equation}
where $\phi_{n,\nu+\mu}$ are the probability densities given by (\ref{eq:sqrtF})
and $c_{n}^{\left(\nu,\mu\right)}$ are positive coefficients that
sum to $1$ defined by 
\begin{equation}
c_{n}^{\left(\nu,\mu\right)}=\frac{1}{B\left(\nu,\mu\right)}\frac{\left(\frac{d}{2}\right)_{n}}{n!}\int_{0}^{1}t^{\mu+\frac{d}{2}-1}\left(1-t\right)^{\nu+\frac{d}{2}-1}\left(1-t+t^{2}\right)^{n}dt,\,\,n\ge0.
\label{eq:cn}
\end{equation}
\end{thm}

Note that (\ref{eq:cn}) shows that $c_{n}^{\left(\mu,\nu\right)} = c_{n}^{\left(\nu,\mu\right)}$ as is to be expected.

\begin{remark}
{\rm The Fourier transform of the Student $t-$density in the case $d=1$ is
\begin{equation}
\hat{f}_{\nu}(y)=\int_{-\infty}^{+\infty} e^{ixy} f_{\nu}\left( x \right) dx = k_{\nu} \left( \vert y \vert \right)
\label{knuy}
\end{equation}
with 
\[
k_{\nu}(u) = \frac{2^{1-\nu}}{\Gamma(\nu)}u^{\nu}K_{\nu}(u),\,\,u\ge0,
\]
where $K_{\nu}$ is the MacDonald function. 

We claim that 
\begin{equation}
D^{2n}\hat{f}_{\mu+\nu+n}\left(y\right)=\left(-1\right)^{n}\frac{\left(\frac{1}{2}\right)_{n}}{\left(\mu+\nu\right)_{n}}\hat{\phi}_{n,\mu+\nu}\left(y\right),\,\,y\in\mathbb{R}.\label{eq:D2n}
\end{equation}
In fact, from (\ref{knuy}), we get
\[
D^{2n}\hat{f}_{\mu+\nu+n}\left(y\right)=\left(-1\right)^{n}\int_{-\infty}^{+\infty}e^{ixy}A_{\mu+\nu+n,1}\frac{x^{2n}}{\left(1+x^{2}\right)^{\mu+\nu+n+\frac{1}{2}}}dx\]
which is (\ref{eq:D2n}) because
\[
A_{\mu+\nu+n,1}B\left(\mu+\nu,n+\frac{1}{2}\right)=\frac{\left(\frac{1}{2}\right)_{n}}{\left(\mu+\nu\right)_{n}}.
\]
Theorem \ref{thm:Student+Student} can be expressed for $d=1$ as
\[
k_{\nu}\left(\vert y\vert\right)k_{\mu}\left(\vert y\vert\right)=\sum_{n=0}^{+\infty}c_{n}^{\left(\nu,\mu\right)}\hat{\phi}_{n,\mu+\nu}\left(y\right)
\]
hence by (\ref{eq:D2n})
\[
=\sum_{n=0}^{+\infty}(-1)^{n}c_{n}^{\left(\nu,\mu\right)}\frac{\left(\mu+\nu\right)_{n}}{\left(\frac{1}{2}\right)_{n}}k_{\mu+\nu+n}^{\left(2n\right)}\left(\vert y\vert\right)
\]
and finally by (\ref{eq:cn}), for $u>0,$ 
\[
B\left(\mu,\nu\right)k_{\nu}\left(u\right)k_{\mu}\left(u\right)=\sum_{n=0}^{+\infty}\frac{\left(-1\right)^{n}\left(\mu+\nu\right)_{n}}{n!}C_{n}(\nu,\mu)k_{\mu+\nu+n}^{\left(2n\right)}\left(u\right)
\]
with 
\[
C_{n}(\nu,\mu) = \int_{0}^{1}t^{\mu-\frac{1}{2}}\left(1-t\right)^{\nu-\frac{1}{2}}\left(1-t+t^{2}\right)^{n}dt.
\]
In terms of the Macdonald function, this equation can be written
\[
K_{\nu}\left(u\right)K_{\mu}\left(u\right)=\sum_{n=0}^{+\infty}\frac{\left(-1\right)^{n}}{n!2^{n+1}}C_{n}(\nu,\mu)u^{-\mu-\nu}D^{2n}\left\{ u^{\mu+\nu+n}K_{\mu+\nu+n}\left(u\right)\right\}. 
\]
We have not been able to find this expression in the literature.\\
Still another interpretation of formula (\ref{fT1T2}) for $d=1$ is to rewrite it as 
\[
f_{\mu}*f_{\nu}\left(x\right)=f_{\mu+\nu}\left(x\right)\sum_{n=0}^{+\infty}\frac{\left(\mu+\nu+\frac{1}{2}\right)_{n}}{\left(\frac{1}{2}\right)_{n}}c_{n}^{\left(\nu,\mu\right)}\left(\frac{x^{2}}{1+x^{2}}\right)^{n},
\]
which shows that $\frac{f_{\mu}*f_{\nu}\left(x\right)}{f_{\mu+\nu}\left(x\right)}$
has a holomorphic extension to the open set $\left\{ z\in\mathbb{C}\,\,\vert\,\,\vert\frac{z^{2}}{1+z^{2}}\vert<1\right\} ,$
which is the region between the two branches of the equilateral hyperbola $y^{2}-x^{2}=\frac{1}{2}.$
}
\end{remark}

\begin{remark}
{\rm In \cite{Ruben}, Ruben shows that in the $d=1$ dimensional case, the density $f_{Z}$ of the convolution $Z=T_{1}+T_{2}$  can be written as a mixture of scaled Student t-distributions with $\mu + \nu $ degrees of freedom, i.e. as
\[
f_{Z}(z)=\int_{2}^{+\infty} \frac{1}{t} f_{\mu + \nu}\left(\frac{z}{t} \right) h(t)dt
\]
for some probability density $h$ on $] 2,+ \infty [$. This representation differs from (\ref{fT1T2}) which is a discrete mixture of the densities (\ref{eq:sqrtF}). 
In the same paper, Ruben obtains several expressions for the density $f_{\mathbf{Z}}$ in the case $\mu = \nu.$ Ruben's formula (3.4) coincides with formula (\ref{fT1T2}) in the special case $\mu = \nu,\,d=1.$}
\end{remark}

\section{Properties of the mixing distributions}

In this section, we derive several properties of the discrete distributions
(\ref{eq:alphakmultivariate}) and (\ref{eq:cn}).

\subsection{Stochastic interpretation}

Since, as a consequence of (\ref{eq:alphakmultivariate}), the coefficients
$\alpha_{k}^{\left(\nu,\gamma\right)}$ are positive, and since they
sum to 1 (see Theorem \ref{thm:thm2}), a stochastic interpretation
of Theorem \ref{thm:Student+Gaussian} is as follows:
\[
\mathbf{Z}=\mathbf{X}+\mathbf{T}\overset{d}{=}\mathbf{G}_{K}
\]
where, given $K=k,$ $\mathbf{G}_{K}$ has the conditional density (\ref{eq:gk})
and $K$ is a discrete random variable with 
\begin{equation}
\Pr\left\{ K=k\right\} =\alpha_{k}^{\left(\nu,\gamma\right)},\,\,\, k=0,1,\dots
\label{eq:K}
\end{equation}

We deduce the following
\begin{prop}
\label{pro:Student+Gaussian}
The sum $\mathbf{Z}=\mathbf{X}+\mathbf{T}$ is distributed as $\sigma\Vert\mathbf{X}_{2K+d}\Vert \mathbf{U}_{d}$
where 
\begin{itemize}
\item vector $\mathbf{U}_d$ is uniformly distributed on the unit sphere $\mathcal{S}_d$ in $\mathbb{R}^d,$
\item vector $\mathbf{X}_{2k+d}\in\mathbb{R}^{2k+d}$ is Gaussian with identity
covariance matrix,
\item random variable $K$ defined by (\ref{eq:K}) follows a compound negative
binomial distribution $NB\left(\frac{d}{2},\frac{a}{1+a}\right),$
where the parameter $a$ is Gamma distributed with scale parameter $\gamma^{2}$
and shape parameter $\nu,$ (in short $a\sim \Gamma\left( \nu,\gamma^2 \right)$) so $a$ has the density
\[
\frac{\gamma^{2\nu}}{\Gamma\left( \nu \right)}e^{-t\gamma^2}t^{\nu - 1},\,\,t>0.
\]
\end{itemize}
\end{prop}
An analogous result in the case of the sum of two Student t-vectors
is as follows: $\mathbf{Y}=\mathbf{T}_{1}+\mathbf{T}_{2}$ is distributed as a random variable
$\Phi_{N}$ such that, given $N=n,$ $\Phi_{N}$ has the
conditional density given by (\ref{eq:sqrtF}) and $N$ is a discrete
random variable defined by 
\begin{equation}
\Pr\left\{ N=n\right\} =c_{n}^{\left(\nu,\mu\right)},\,\,\,n=0,1,\dots\label{eq:N}\end{equation}
We deduce the following
\begin{prop}
\label{pro:Gaussian+Gaussian}The sum $\mathbf{Y}=\mathbf{T}_{1}+\mathbf{T}_{2}$ is distributed
as $\Vert\mathbf{T}_{2N+d}\Vert \mathbf{U}_{d}$ where
\begin{itemize}
\item vector $\mathbf{U}_d$ is uniformly distributed on the sphere $\mathcal{S}_d$ in $\mathbb{R}^d,$
\item random vector $\mathbf{T}_{2n+d}\in\mathbb{R}^{2n+d}$ is Student
t-distributed with $2\nu+2\mu$ degrees of freedom
\item random variable $N$ defined by (\ref{eq:N}) follows a compound negative
binomial distribution $NB\left(\frac{d}{2},t\left(1-t\right)\right),$
where the parameter $t$ follows a Beta distribution with parameters $\mu$
and $\nu.$
\end{itemize}
\end{prop}
\subsection{Hausdorff moment sequence}

A sequence $\left\{ s_{n}\right\} $ is called a Hausdorff moment sequence
if it has the form
\[
s_{n}=\int_{0}^{1}x^{n}d\mu\left(x\right),\,\, n=0,1,\dots
\]
for some positive measure $\mu.$ It is called normalized if $\mu_{0}=1.$ The following result holds:
\begin{thm}
\label{thm:Hausdorff}
Both sequences $\left\{ \alpha_{k}^{\left(\nu,\gamma\right)}\right\} $
and $\left\{ c_{n}^{\left(\nu,\mu\right)}\right\} $ are Hausdorff
moment sequences for $d=1$ and $d=2.$
\end{thm}

\subsection{Mean and Variance}

The first and second order moments of these distributions can be easily
computed as a consequence of the fact that they are compound distributions.
\begin{thm}
\label{thm:thm2}
The distribution (\ref{eq:K}) of $K$ verifies
\begin{eqnarray*}
\sum_{k=0}^{+\infty}\alpha_{k}^{\left(\nu,\gamma\right)} & = & 1,\\
EK & = & \frac{d\gamma^{2}}{2\nu-2},\,\,\nu>1,\\
\sigma_{K}^{2} & = & \frac{d\gamma^{2}}{2\left(\nu-1\right)}\left(1+\frac{\gamma^{2}}{2}\frac{d+2\left( \nu-1 \right)}{\left(\nu-1\right)\left(\nu-2\right)}\right),\,\,\nu>2.
\end{eqnarray*}
\end{thm}
In the case of the sum of Student t- variables, we have
\begin{thm}
\label{thm:thm2-1}
The distribution (\ref{eq:N}) of $N$ verifies\begin{eqnarray*}
\sum_{k=0}^{+\infty}c_{n}^{\left(\nu,\mu\right)} & = & 1,\\
EN & = & \frac{d}{2}\left(\frac{B\left(\nu-1,\mu-1\right)}{B\left(\nu,\mu\right)}-1\right)\\
\sigma_{N}^{2} & = & \frac{d}{4}\left(\left( d+2 \right)\frac{B\left(\nu-2,\mu-2\right)}{B\left(\nu,\mu\right)}-2\frac{B\left(\nu-1,\mu-1\right)}{B\left(\nu,\mu\right)}-d\frac{B^{2}\left(\nu-1,\mu-1\right)}{B^{2}\left(\nu,\mu\right)}\right).\end{eqnarray*}
\end{thm}

\subsection{Asymptotic behaviors}

\subsubsection{fixed value of $\gamma$}

Let us now study the asymptotic behavior of the distribution (\ref{eq:fzmultivariate})
for large value of the degrees of freedom $\nu:$ note that the two
next theorems can be checked directly by looking at the stochastic
representation of $K$ as given in Proposition \ref{pro:Student+Gaussian},
but their proofs are instructive by themselves.
\begin{thm}
\label{thm:thm3}
For a fixed value of the parameter $\gamma,$
\[
\lim_{\nu\to+\infty}\alpha_{k}^{\left(\nu,\gamma\right)}=\begin{cases}
1 & k=0\\
0 & k\ne0\end{cases}
\]
\end{thm}
We note that this result is a consequence of the fact that since $\mathbf{T}$
has covariance matrix $\frac{1}{2\nu-2}I_{d},$ the random variable $\mathbf{Z}=\mathbf{X}+\mathbf{T}$ has the
same distribution as $\mathbf{X}$ as $\nu\to+\infty.$

\subsubsection{covariance matrix of $\mathbf{T}$ fixed}

Let us now look at the asymptotic behavior of these coefficients for
large values of $\nu$ but for a fixed value of the covariance matrix of $\mathbf{T}:$
since $\mathbf{T}$ has covariance matrix $\frac{1}{2\left(\nu-1\right)}I_{d},$ let us consider
\begin{equation}
\mathbf{Z}=\mathbf{X}+\eta\sqrt{2\left(\nu-1\right)}\mathbf{T},
\label{eq:Zfixedvariance}
\end{equation}
where $\mathbf{X}$ has covariance matrix $\sigma^{2}I_{d}$ and $\eta\sqrt{2\left(\nu-1\right)}\mathbf{T}$
has covariance matrix $\eta$$^{2}I_{d}$. Define
\[
\tilde{\gamma}=\frac{\eta\sqrt{\nu-1}}{\sigma}
\]
so that 
\begin{equation}
f_{\mathbf{Z}}\left(\mathbf{z}\right)=\sum_{k=0}^{+\infty}\alpha_{k}^{\left(\nu,\tilde{\gamma}\right)}g_{k,\sigma}\left(\mathbf{z}\right).
\label{eq:fZfixedvariance}
\end{equation}

\begin{thm}
\label{thm:thm4}
The coefficients $\alpha_{k}^{\left(\nu,\tilde{\gamma}\right)}$
in (\ref{eq:fZfixedvariance}) behave for large $\nu$ as
\[
\alpha_{k}^{\left(\nu,\tilde{\gamma}\right)}\sim\frac{\Gamma\left(k+\frac{d}{2}\right)}{k!\Gamma\left( \frac{d}{2} \right)}\left(\frac{\eta^{2}}{\sigma^{2}}\right)^{k}\left(1+\frac{\eta^{2}}{\sigma^{2}}\right)^{-k-\frac{d}{2}}
\]
As a consequence, the density of $\mathbf{Z}=\mathbf{X}+\eta\sqrt{2\left(\nu-1\right)}\mathbf{T}$ converges to a Gaussian
density with covariance matrix $\left( \sigma^{2}+\eta^{2} \right)I_{d}$, as expected.
\end{thm}

\begin{figure}[h] 
   \centering
   \includegraphics[scale=0.3]{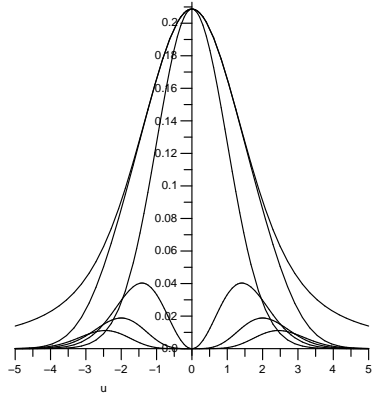} 
   \caption{from top to bottom, density $f_{Z}$ in the case $d=1$, its approximation by the sum
of the first four term in (\ref{eq:fzmultivariate}) and these first
four terms}
   \label{Flo:fig1}
\end{figure}

In Fig.\ref{Flo:fig1} are shown the succesive terms of the approximation
of the density of $Z$ in the one-dimensional case ($d=1$): the upper curve is the density $f_{Z},$ the
second one is the approximation of $f_{Z}$ by the sum of the four
first terms $0\le k\le3$ in (\ref{eq:fzmultivariate}), and the four
lower curves are the successive terms of this sum for $k=0,1,2$ and
$3.$ The first $(k=0)$ term (Gaussian density) fits $f_{Z}$ around
$z=0,$ but the tail is underestimated; the other terms contribute
to correct this underestimation, each adding, as $k$ increases, a
smaller contribution farther from the origin.

\section{\label{sec:Proofs}Proofs}

\subsection{Proof of Theorem \ref{thm:Student+Gaussian}}

Using a result by Keilson and Steutel \cite{Keilson}, a stochastic
representation for the random variable $\mathbf{Z}$ is 
\[
\mathbf{Z}=\sigma \mathbf{N}_{1}+\frac{1}{\sqrt{2a}}\mathbf{N}_{2}\overset{d}{=}\sqrt{\sigma^{2}+\frac{1}{2a}}\mathbf{N}
\]
where $\mathbf{N}$ is $d-$variate standard normal, $a$ is Gamma distributed $\Gamma\left(\nu,1\right)$
and $\overset{d}{=}$ denotes equality in distribution. The density
of $\mathbf{Z}$ reads thus
\begin{equation}
f_{\mathbf{Z}}\left(\mathbf{z}\right)=\int_{0}^{+\infty}\frac{1}{\Gamma\left(\nu\right)}e^{-a}a^{\nu-1}\frac{1}{\left( 2\pi\left(\sigma^{2}+\frac{1}{2a}\right) \right)^{\frac{d}{2}}}\exp\left(-\frac{\Vert\mathbf{z}\Vert^{2}}{2\left(\sigma^{2}+\frac{1}{2a}\right)}\right)da.
\end{equation}
Since
$\exp\left(-\frac{\Vert\mathbf{z}\Vert^{2}}{2\left(\sigma^{2}+\frac{1}{2a}\right)}\right)=\exp\left(-\frac{\Vert\mathbf{z}\Vert^{2}}{2\sigma^{2}}\right)\exp\left(\frac{\Vert\mathbf{z}\Vert^{2}}{2\sigma^{2}}\frac{\gamma^{2}}{a+\gamma^{2}}\right)$
with $\gamma^{2}=\frac{1}{2\sigma^{2}}$, we deduce
\begin{eqnarray*}
f_{\mathbf{Z}}\left(\mathbf{z}\right) & = & \frac{1}{\left( \sigma\sqrt{2\pi} \right)^{d}}\exp\left(-\frac{\Vert \mathbf{z}\Vert ^{2}}{2\sigma^{2}}\right)\frac{1}{\Gamma\left(\nu\right)}\int_{0}^{+\infty}e^{-a}a^{\nu-1}\frac{a^{\frac{d}{2}}}{\left( a+\gamma^{2} \right)^{\frac{d}{2}}}\exp\left(\frac{\Vert \mathbf{z}^{2}\Vert}{2\sigma^{2}}\frac{\gamma^{2}}{a+\gamma^{2}}\right)da\\
 & = & \frac{1}{\left( \sigma\sqrt{2\pi} \right)^{d}}\exp\left(-\frac{\Vert \mathbf{z}\Vert^{2}}{2\sigma^{2}}\right)\frac{1}{\Gamma\left(\nu\right)}\sum_{k=0}^{+\infty}\frac{\gamma^{2k}}{k!}\left(\frac{\Vert \mathbf{z}\Vert^{2}}{2\sigma^{2}}\right)^{k}\int_{0}^{+\infty}e^{-a}a^{\nu+\frac{d}{2}-1}\left(a+\gamma^{2}\right)^{-k-\frac{d}{2}}da,
\end{eqnarray*}
where the interchange of the sum and the integral is justified by the
Lebesgue monotone convergence theorem. Since, using \cite[4.642]{Gradshteyn}, $\forall k\in\mathbb{N},$
\[
\int_{\mathbb{R}^{d}}\left(\frac{\Vert\mathbf{x}\Vert^{2}}{2\sigma^{2}}\right)^{k}\exp\left(-\frac{\Vert\mathbf{x}\Vert^{2}}{2\sigma^{2}}\right)d\mathbf{x}=\frac{2\pi^{\frac{d}{2}}}{\Gamma\left(\frac{d}{2}\right)}\int_{0}^{+\infty}\left(\frac{r^{2}}{2\sigma^{2}}\right)^{k}\exp\left(-\frac{r^{2}}{2\sigma^{2}}\right)r^{d-1}dr,
\]
which is equal to $\frac{\pi^{\frac{d}{2}}}{\Gamma\left(\frac{d}{2}\right)}\left(\sigma\sqrt{2}\right)^{d}\Gamma\left(k+\frac{d}{2}\right),$
the function $g_{k,\sigma}$ as defined by (\ref{eq:gk}) is a probability
density and the result follows. 
Using the Kummer function 
\begin{equation}
\Psi(\alpha,\beta;z)=\frac{1}{\Gamma(\alpha)}\int_{0}^{\infty} e^{-zt} t^{\alpha-1}(1+t)^{\beta - \alpha -1} dt,
\end{equation}
cf. \cite[p.1085]{Gradshteyn}, the integral in the expression of $f_{\mathbf{Z}}$ can be written 
\begin{equation}
\int_{0}^{\infty}a^{\nu + \frac{d}{2}-1} e^{-a}(a+\gamma^2)^{-k-\frac{d}{2}}da = (\gamma^2)^{\nu -k}\Gamma(\nu + \frac{d}{2})\Psi(\nu + \frac{d}{2},\nu -k+1;\gamma^2).
\end{equation}

\subsection{Proof of Proposition \ref{pro:extensionStudentGaussian}}

We have
\[
\tilde{\mathbf{Z}}=a_{2}\left(\frac{a_{1}}{a_{2}}\mathbf{X}+\mathbf{T}\right)=a_{2}\left(\tilde{\sigma}\mathbf{N}_{1}+\mathbf{T}\right)=a_{2}\tilde{\mathbf{Z}}
\]
with $\tilde{\sigma}=\frac{ a_{1}}{ a_{2}}\sigma$
so that $\tilde{\gamma}=\frac{ a_{2}}{\sigma\sqrt{2} a_{1}}$
and 
\[
f_{\tilde{\mathbf{Z}}}\left(z\right)=\frac{1}{a_{2}^{d}}f_{\mathbf{Z}}\left(\frac{\mathbf{z}}{a_{2}}\right)=\sum_{k=0}^{+\infty}\alpha_{k}^{\left(\nu,\tilde{\gamma}\right)}\frac{1}{a_{2}^{d}}g_{k,\frac{a_{1}}{a_{2}}\sigma}\left(\frac{\mathbf{z}}{a_{2}}\right)
\]
but $\frac{1}{a_{2}^{d}}g_{k,\frac{a_{1}}{a_{2}}\sigma}\left(\frac{\mathbf{z}}{a_{2}}\right)=g_{k,a_{1}\sigma}\left(\mathbf{z}\right)$
so that the result holds.

\subsection{Proof of Theorem \ref{thm:Student+Student}}

Let us define $t=\frac{\sigma^{2}}{2}$ so that $\gamma=\frac{1}{\sigma\sqrt{2}}=\frac{1}{2\sqrt{t}}$
and consider the $k-$th term in the convex combination (\ref{eq:fzmultivariate})
\[
\alpha_{k}^{\left(\nu,\frac{1}{2\sqrt{t}}\right)}g_{k,\sqrt{2t}}\left(\mathbf{x}\right)=\frac{1}{\pi^{\frac{d}{2}}\Gamma\left(\nu\right)k!}\left(\frac{\Vert \mathbf{x}\Vert ^{2}}{4t}\right)^{k}\exp\left(-\frac{\Vert \mathbf{x}\Vert ^{2}}{4t}\right)\int_{0}^{+\infty}\frac{s^{\nu+\frac{d}{2}-1}\exp\left(-s\right)}{\left(1+4ts\right)^{k+\frac{d}{2}}}ds;
\]
putting $w=4ts$ and subordinating with the inverse Gamma distribution
yields
\begin{eqnarray*}
\int_{0}^{+\infty}\alpha_{k}^{\left(\nu,\frac{1}{2\sqrt{t}}\right)}g_{k,\sqrt{2t}}\left(\mathbf{x}\right)d\tilde{\gamma}_{\mu}\left(t\right) & = & \frac{2^{-2\mu}}{\pi^{\frac{d}{2}}\Gamma\left(\mu\right)\Gamma\left(\nu\right)}\frac{\Vert \mathbf{x} \Vert^{2k}}{k!}
\int_{0}^{+\infty}\left(\frac{1}{4t}\right)^{k+\nu+\frac{d}{2}}\frac{e^{-\frac{\Vert \mathbf{x} \Vert^{2}+1}{4t}}}{t^{\mu+1}}\int_{0}^{+\infty}\frac{w^{\nu+\frac{d}{2}-1}e^{-\frac{w}{4t}}}{\left(1+w\right)^{k+\frac{d}{2}}}dwdt\\
 & = & \frac{4^{-\mu-\nu-k-\frac{d}{2}}}{\pi^{\frac{d}{2}}\Gamma\left(\mu\right)\Gamma\left(\nu\right)}\frac{\Vert \mathbf{x} \Vert^{2k}}{k!}\int_{0}^{+\infty}\frac{w^{\nu+\frac{d}{2}-1}}{\left(1+w\right)^{k+\frac{d}{2}}}\int_{0}^{+\infty}\frac{e^{-\frac{w+\Vert \mathbf{x} \Vert^{2}+1}{4t}}}{t^{k+\mu+\nu+\frac{d}{2}+1}}dtdw.
 \end{eqnarray*}
By the change of variable $u=\frac{w+\Vert \mathbf{x}\Vert^{2}+1}{4t},$ this is equal
to 
\[
\frac{1}{4^{\mu+\nu+k+\frac{d}{2}}\pi^{\frac{d}{2}}\Gamma\left(\mu\right)\Gamma\left(\nu\right)}\frac{\Vert \mathbf{x}\Vert ^{2k}}{k!}\int_{0}^{+\infty}\frac{w^{\nu+\frac{d}{2}-1}}{\left(1+w\right)^{k+\frac{d}{2}}}\left(\frac{w+\Vert \mathbf{x}\Vert^{2}+1}{4}\right)^{-k-\mu-\nu-\frac{d}{2}}\int_{0}^{+\infty}u^{k+\mu+\nu+\frac{d}{2}-1}e^{-u}dudw
\]
and, since the latter integral is a Gamma integral, to
\[
\frac{\Gamma\left(\mu+\nu+k+\frac{d}{2}\right)}{\pi^{\frac{d}{2}}\Gamma\left(\mu\right)\Gamma\left(\nu\right)}\frac{\Vert \mathbf{x} \Vert^{2k}}{k!}\int_{0}^{+\infty}\frac{w^{\nu+\frac{d}{2}-1}}{\left(1+w\right)^{k+\frac{d}{2}}\left(w+\Vert \mathbf{x} \Vert^{2}+1\right)^{k+\mu+\nu+\frac{d}{2}}}dw
\]
The latter integral can be identified as Euler's integral representation
of a hypergeometric function \cite[Th. 2.2.1 p.65]{Andrews} 
\begin{eqnarray*}
 & \int_{0}^{+\infty}\frac{w^{\nu+\frac{d}{2}-1}}{\left(1+w\right)^{k+\frac{d}{2}}\left(w+\Vert \mathbf{x} \Vert^{2}+1\right)^{k+\mu+\nu+\frac{d}{2}}}dw
\\
= & \frac{\Gamma\left(2k+\mu+\frac{d}{2}\right)\Gamma\left(\nu+\frac{d}{2}\right)}{\Gamma\left(2k+\mu+\nu+d\right)}{}_{2}F_{1}\left(\begin{array}{c}
k+\mu+\nu+\frac{d}{2},\,\,2k+\mu+\frac{d}{2}\\
2k+\mu+\nu+d\end{array};-\Vert \mathbf{x}\Vert^{2}\right)
\end{eqnarray*}
so that, using Pfaff's formula \cite[2.2.6 p.68]{Andrews} 
\[
_{2}F_{1}\left(\begin{array}{c}
a,\,\,b\\
c\end{array};x\right)=\left(1-x\right)^{-a}{}_{2}F_{1}\left(\begin{array}{c}
a,\,c-b\\
c\end{array};\frac{x}{x-1}\right)
\]
we deduce
\begin{eqnarray*}
\int_{0}^{+\infty}\alpha_{k}^{\left(\nu,\frac{1}{2\sqrt{t}}\right)}g_{k,\sqrt{2t}}\left(\mathbf{x}\right)d\tilde{\gamma}_{\mu}\left(t\right) & = & \frac{\Gamma\left(\mu+\nu+k+\frac{d}{2}\right)\Gamma\left(\nu+\frac{d}{2}\right)}{\pi^{\frac{d}{2}}\Gamma\left(\mu\right)\Gamma\left(\nu\right)k!}\frac{\Gamma\left(2k+\mu+\frac{d}{2}\right)}{\Gamma\left(2k+\mu+\nu+d\right)}\\
 & \times & \frac{\Vert \mathbf{x} \Vert^{2k}}{\left(1+\Vert \mathbf{x} \Vert^{2}\right)^{k+\mu+\nu+\frac{d}{2}}}{}_{2}F_{1}\left(\begin{array}{c}
k+\mu+\nu+\frac{d}{2},\,\,\nu+\frac{d}{2}\\
2k+\mu+\nu+d\end{array};\frac{\Vert \mathbf{x}\Vert^{2}}{1+\Vert \mathbf{x} \Vert^{2}}\right).
\end{eqnarray*}
Replacing the hypergeometric function by its series, we deduce
\begin{eqnarray*}
\int_{0}^{+\infty}\alpha_{k}^{\left(\nu,\frac{1}{2\sqrt{t}}\right)}g_{k,\sqrt{2t}}\left(\mathbf{x}\right)d\tilde{\gamma}_{\mu}\left(t\right) & = & \frac{\Gamma\left(\mu+\nu+k+\frac{d}{2}\right)\Gamma\left(\nu+\frac{d}{2}\right)}{\pi^{\frac{d}{2}}\Gamma\left(\mu\right)\Gamma\left(\nu\right)k!}\frac{\Gamma\left(2k+\mu+\frac{d}{2}\right)}{\Gamma\left(2k+\mu+\nu+d\right)}\\
 & \times & \sum_{l=0}^{+\infty}\frac{\left(\nu+\frac{d}{2}\right)_{l}\left(k+\mu+\nu+\frac{d}{2}\right)_{l}}{\left(2k+\mu+\nu+d\right)_{l}l!}\frac{\Vert \mathbf{x} \Vert^{2k+2l}}{\left(1+\Vert \mathbf{x} \Vert^{2}\right)^{k+l+\mu+\nu+\frac{d}{2}}}\\
 & = & \sum_{l=0}^{+\infty}d_{k,l}\phi_{l+k,\eta}\left(\mathbf{x}\right)
\end{eqnarray*}
where $\eta=\mu+\nu,$ the density $\phi_{n,\eta}\left(\mathbf{x}\right)$ is defined
by (\ref{eq:sqrtF}) and the coefficient $d_{k,l}$ is equal to
\begin{eqnarray*}
d_{k,l} & = & \frac{\Gamma\left(\mu+\nu+k+\frac{d}{2}\right)\Gamma\left(\nu+\frac{d}{2}\right)}{\pi^{\frac{d}{2}}\Gamma\left(\mu\right)\Gamma\left(\nu\right)k!}\frac{\Gamma\left(2k+\mu+\frac{d}{2}\right)}{\Gamma\left(2k+\mu+\nu+d\right)}\frac{\left(\nu+\frac{d}{2}\right)_{l}\left(k+\mu+\nu+\frac{d}{2}\right)_{l}}{\left(2k+\mu+\nu+d\right)_{l}l!}\frac{\Gamma\left(\mu+\nu\right)\Gamma\left(k+l+\frac{d}{2}\right)}{\Gamma\left(\mu+\nu+k+l+\frac{d}{2}\right)}\frac{\pi^{\frac{d}{2}}}{\Gamma\left( \frac{d}{2} \right)}\\
 & = & \frac{\Gamma\left(\nu+\frac{d}{2}\right)\left(\nu+\frac{d}{2}\right)_{l}}{\Gamma\left(\mu\right)\Gamma\left(\nu\right)k!l!}\frac{\Gamma\left(2k+\mu+\frac{d}{2}\right)\Gamma\left(\mu+\nu\right)\Gamma\left(k+l+\frac{d}{2}\right)}{\Gamma\left(2k+l+\mu+\nu+d\right)\Gamma\left(\frac{d}{2}\right)}.
\end{eqnarray*}
using the identity $\left(a\right)_{l}\Gamma\left(a\right)=\Gamma\left(a+l\right)$
twice. Writing $\Gamma\left(2k+\mu+\frac{d}{2}\right)=\left(\mu+\frac{d}{2}\right)_{2k}\Gamma\left(\mu+\frac{d}{2}\right),$
$\Gamma\left(k+l+\frac{d}{2}\right)=\left(\frac{d}{2}\right)_{k+l}\Gamma\left(\frac{d}{2}\right)$
and $\Gamma\left(2k+l+\mu+\nu+d\right)=\left(\mu+\nu+d\right)_{2k+l}\Gamma\left(\mu+\nu+d\right),$
we deduce
\[
d_{k,l}=\frac{B\left(\nu+\frac{d}{2},\mu+\frac{d}{2}\right)}{B\left(\nu,\mu\right)k!l!}\frac{\left(\mu+\frac{d}{2}\right)_{2k}\left(\nu+\frac{d}{2}\right)_{l}\left(\frac{d}{2}\right)_{k+l}}{\left(\mu+\nu+d\right)_{2k+l}}.
\]
Now the density of $\mathbf{Y}=\mathbf{T}_{1}+\mathbf{T}_{2}$ is the sum over $k$ of these
terms and reads
\[
f_{Y}\left(\mathbf{x}\right)=\sum_{k,l=0}^{+\infty}d_{k,l}\phi_{k+l,\eta}\left(\mathbf{x}\right).\]
We perform the change of variable $n=k+l$ to obtain
\[
f_{\mathbf{Y}}\left(\mathbf{x}\right)=\sum_{n=0}^{+\infty}\phi_{n,\eta}\left(\mathbf{x}\right)\sum_{k=0}^{n}d_{k,n-k}.
\]
The inner sum can be written
\begin{eqnarray*}
c_{n}^{\left(\nu,\mu\right)}=\sum_{k=0}^{n}d_{k,n-k} & = & \frac{B\left(\nu+\frac{d}{2},\mu+\frac{d}{2}\right)}{B\left(\nu,\mu\right)}\frac{\left(\frac{d}{2}\right)_{n}}{n!}\sum_{k=0}^{n}\binom{n}{k}\frac{\left(\mu+\frac{d}{2}\right)_{2k}\left(\nu+\frac{d}{2}\right)_{n-k}}{\left(\mu+\nu+d\right)_{n+k}}\\
  =  \frac{B\left(\nu+\frac{d}{2},\mu+\frac{d}{2}\right)}{B\left(\nu,\mu\right)}&\frac{\left(\frac{d}{2}\right)_{n}}{\left(\mu+\nu+d\right)_{2n}n!}&\sum_{k=0}^{n}\binom{n}{k}\left(\mu+\frac{d}{2}\right)_{2k}\left(\nu+\frac{d}{2}\right)_{n-k}\left(\mu+\nu+d+n+k\right)_{n-k}
\end{eqnarray*}
where we have used the identity
\[
\frac{1}{\left(\mu+\nu+d\right)_{n+k}}=\frac{\left(\mu+\nu+d+n+k\right)_{n-k}}{\left(\mu+\nu+d\right)_{2n}}.
\]
By the Chu-Vandermonde identity \cite[p.70]{Andrews} 
\begin{eqnarray*}
\left(\mu+\nu+d+n+k\right)_{n-k} & = & \left(\left(\mu+\frac{d}{2}+2k\right)+\left(\nu+\frac{d}{2}+n-k\right)\right)_{n-k}\\
 & = & \sum_{j=0}^{n-k}\binom{n-k}{j}\left(\mu+\frac{d}{2}+2k\right)_{n-k-j}\left(\nu+\frac{d}{2}+n-k\right)_{j}
 \end{eqnarray*}
and using the identities $\left(\mu+\frac{d}{2}\right)_{2k}\left(\mu+\frac{d}{2}+2k\right)_{n-k-j}=\left(\mu+\frac{d}{2}\right)_{n+k-j}$
and $\left(\nu+\frac{d}{2}\right)_{n-k}\left(\nu+\frac{d}{2}+n-k\right)_{j}=\left(\nu+\frac{d}{2}\right)_{n-k+j},$
we deduce
\begin{eqnarray*}
c_{n}^{\left(\nu,\mu\right)} & = & \frac{B\left(\nu+\frac{d}{2},\mu+\frac{d}{2}\right)}{B\left(\nu,\mu\right)}\frac{\left(\frac{d}{2}\right)_{n}}{\left(\mu+\nu+d\right)_{2n}n!}\sum_{k=0}^{n}\sum_{j=0}^{n-k}\binom{n}{k}\binom{n-k}{j}\left(\mu+\frac{d}{2}\right)_{n+k-j}\left(\nu+\frac{d}{2}\right)_{n-k+j}\\
 & = & \frac{B\left(\nu+\frac{d}{2},\mu+\frac{d}{2}\right)}{B\left(\nu,\mu\right)}\frac{\left(\frac{d}{2}\right)_{n}}{n!}\sum_{p=-n}^{n}\binom{n}{p}_{2}\frac{\left(\mu+\frac{d}{2}\right)_{n+p}\left(\nu+\frac{d}{2}\right)_{n-p}}{\left(\mu+\nu+d\right)_{2n}},
\end{eqnarray*}
where we have used the change of summation index $p=k-j$ and where
the trinomial \cite{Andrews2} coefficient $\binom{n}{p}_{2}$ is\[
\binom{n}{p}_{2}=\sum_{\underset{j+k\le n;\,k-j=p}{j,k=0}}^{n}\binom{n}{k}\binom{n-k}{j},\,\,-n\le p\le n.\]

Writing $\left(\mu+\frac{d}{2}\right)_{n+p}\Gamma\left(\mu+\frac{d}{2}\right)=\Gamma\left(\mu+n+p+\frac{d}{2}\right)$
and $\left(\nu+\frac{d}{2}\right)_{n-p}\Gamma\left(\nu+\frac{d}{2}\right)=\Gamma\left(\nu+\frac{d}{2}+n-p\right)$
and replacing the resulting Beta function by its integral representation
we obtain
\[
c_{n}^{\left(\nu,\mu\right)}=\frac{1}{B\left(\nu,\mu\right)}\frac{\left(\frac{d}{2}\right)_{n}}{n!}\sum_{p=-n}^{n}\binom{n}{p}_{2}\int_{0}^{1}t^{\mu+n+p+\frac{d}{2}-1}\left(1-t\right)^{\nu+n-p+\frac{d}{2}-1}dt,
\]
and since the generating function of the trinomial coefficient is
\[
\sum_{p=-n}^{+n}\binom{n}{p}_{2}x^{p}=\left(x+x^{-1}+1\right)^{n}
\]
we finally obtain
\[
c_{n}^{\left(\nu,\mu\right)}=\frac{1}{B\left(\nu,\mu\right)}\frac{\left(\frac{d}{2}\right)_{n}}{n!}\int_{0}^{1}t^{\mu+n+\frac{d}{2}-1}\left(1-t\right)^{\nu+n+\frac{d}{2}-1}\left(1+\frac{t}{1-t}+\frac{1-t}{t}\right)^{n}dt
\]
which is the desired result.

\subsection{Proof of Theorem \ref{thm:Hausdorff}}

In the case of the coefficients $\alpha_{k}^{\left(\nu,\gamma\right)},$
we can factorize
\[
\alpha_{k}^{\left(\nu,\gamma\right)}=\frac{\left( \frac{d}{2} \right)_{k}}{k!}\tau_{k}
\]
where 
\[
\tau_{k}=\frac{\left( \sigma\sqrt{2} \right)^d}{\Gamma\left(\nu\right)}\int_{0}^{+\infty}\frac{s^{\nu+\frac{d}{2}-1}e^{-s}}{\left(1+2\sigma^{2}s\right)^{k+\frac{d}{2}}}ds=\frac{1}{\left(2\sigma^{2}\right)^{\nu}\Gamma\left(\nu\right)}\int_{0}^{1}u^{k}\frac{\left(1-u\right)^{\nu+\frac{d}{2}-1}}{u^{\nu+1}}\exp\left(-\frac{1-u}{2\sigma^{2}u}\right)du.
\]

This expression shows that the sequence $\left\{ \tau_{k}\right\} $
is also a Hausdorff moment sequence (not normalized) and hence $\alpha_{k}^{\left(\nu,\gamma\right)}$
is a Hausdorff moment sequence for $d=1,2$ as product of two such sequences. It
is not normalized since $\sum_{k=0}^{+\infty}\alpha_{k}^{\left(\nu,\gamma\right)}=1.$
In fact, $\frac{\left( \frac{d}{2} \right)_{k}}{k!}$ is trivially a Hausdorff moment sequence when $d=2,$ and for $d=1$ we have 
\[
\frac{\left( \frac{1}{2} \right)_{k}}{k!}=\frac{1}{\pi}\int_{0}^{1}\frac{s^{k}ds}{\sqrt{s\left(1-s\right)}}.
\]
It is a general fact that $s_{k}=\frac{\left( a \right)_{k}}{\left( b
  \right)_{k}}$ is a Hausdorff moment sequence when $0<a\le b,$ see
\cite[formula (26)]{Berg} and if $a>b$ then $s_{k}$ is not even a Hamburger moment sequence since the Hankel determinant $D_{2}=s_{0}s_{2}-s_{1}^{2}=\frac{a\left( b-a \right)}{b^{2}\left( b+1 \right)}<0.$

This shows that the density (\ref{eq:fzmultivariate}) is an infi{}nite
convex combination of the densities (\ref{eq:gk}) and the discrete probability
$\alpha^{\left(\nu,\gamma\right)}=\sum_{k=0}^{+\infty}\alpha_{k}^{\left(\nu,\gamma\right)}\delta_{k}$
is infinitely divisible for $d=1,2$ since the sequence $\left\{ \alpha_{k}^{\left(\nu,\gamma\right)}\right\} $
is completely monotonic, which is the same as being a Hausdorff moment
sequence.


In the case of the sequence $\left\{ c_{n}^{\left(\nu,\mu\right)}\right\} ,$
we find
\begin{eqnarray*}
c_{n}^{\left(\nu,\mu\right)} & = &
\frac{1}{B\left(\nu,\mu\right)}\frac{\left(\frac{d}{2}\right)_{n}}{n!}\int_{\frac{3}{4}}^{1}u^{n}
\times\\
 &  & \left[\left(\frac{1}{2}-\sqrt{u-\frac{3}{4}}\right)^{\mu+\frac{d}{2}-1}\left(\frac{1}{2}+\sqrt{u-\frac{3}{4}}\right)^{\nu+\frac{d}{2}-1}+\left(\frac{1}{2}+\sqrt{u-\frac{3}{4}}\right)^{\mu+\frac{d}{2}-1}\left(\frac{1}{2}-\sqrt{u-\frac{3}{4}}\right)^{\nu+\frac{d}{2}-1}\right]\frac{du}{2\sqrt{u-\frac{3}{4}}}.
\end{eqnarray*}
As before
 $\left\{ c_{n}^{\left(\nu,\mu\right)}\right\} $
is a Hausdorff moment sequence for $d=1,2$.

\subsection{Proof of Theorem \ref{thm:thm2}}

Although the first result can be obtained by integrating (\ref{eq:fzmultivariate})
over $\mathbf{z}\in\mathbb{R}^{d},$ we check it directly since the computation
of the sum is instructive by itself. For any $u\in\left]-1,1\right[,$
we have by the binomial series
\begin{equation}
\sum_{k=0}^{+\infty}\frac{\Gamma\left(k+\frac{d}{2}\right)}{k!}u^{k}=\frac{\Gamma\left( \frac{d}{2} \right)}{\left( 1-u \right)^\frac{d}{2}}
\label{eq:binom}
\end{equation}
so that, with $u=\frac{\gamma^{2}}{a+\gamma^{2}}$,\[
\sum_{k=0}^{+\infty}\frac{\Gamma\left(k+\frac{d}{2}\right)}{k!}\left(\frac{\gamma^{2}}{a+\gamma^{2}}\right)^{k}=\Gamma\left( \frac{d}{2} \right)\left( \frac{a+\gamma^{2}}{a} \right)^\frac{d}{2}\]
and we deduce
\begin{eqnarray*}
\sum_{k=0}^{+\infty}\alpha_{k}^{\left( \nu,\gamma \right)} &=& \int_{0}^{+\infty}\frac{e^{-a}}{\Gamma\left(\nu\right)\Gamma\left( \frac{d}{2} \right)} 
\frac{a^{\nu+\frac{d}{2}-1}}{\left(a+\gamma^{2}\right)^\frac{d}{2}} \Gamma\left( \frac{d}{2} \right)
\left(\frac{a+\gamma^{2}}{a}\right)^\frac{d}{2}da =1.
\end{eqnarray*}

The expectation and variances are easily computed from the stochastic
representation of $K=NB\left(\frac{d}{2},\frac{a}{1+a}\right);$
the mean of a negative binomial variable $NB\left(r,p\right)$ being
$\frac{r\left(1-p\right)}{p},$ 
\[
EK=\frac{d}{2}E\frac{\frac{1}{1+a}}{\frac{a}{1+a}}=\frac{d}{2}E\frac{1}{a}=\frac{d\gamma^{2}}{2\nu-2}
\]
since $a$ is $\Gamma\left(\nu,\gamma^{2}\right).$ Since moreover,
the second order moment of a negative binomial is $\frac{r\left(1-p\right)\left(1+r\left(1-p\right)\right)}{p^{2}},$
a straightforward computation gives
\begin{eqnarray*}
EK^{2} & = & \frac{d}{2}\left( E\frac{1+a}{a^{2}}+\frac{d}{2}E\frac{1}{a^{2}} \right)=\frac{d}{2}\left(1+\frac{d}{2}  \right)E\frac{1}{a^{2}}+\frac{d}{2}E\frac{1}{a}\\
 & = & \frac{d}{2} \left(1+\frac{d}{2}  \right)\frac{\gamma^{4}}{\left(\nu-1\right)\left(\nu-2\right)}+\frac{d}{2}\frac{\gamma^{2}}{\nu-1} 
 \end{eqnarray*}
and 
\[
\sigma_{K}^{2}=\frac{d}{2}\frac{\gamma^{2}}{\nu-1}\left(1+\frac{\gamma^{2}}{2}\frac{d+2\left( \nu-1 \right)}{\left(\nu-1\right)\left(\nu-2\right)}\right).\]

\subsection{Proof of Theorem \ref{thm:thm2-1}}

The fact that the coefficients $c_{n}^{\left(\nu,\mu\right)}$ sum to
$1$ is proved in the same way as in the Proof of Theorem \ref{thm:thm2}; the expectation is computed accordingly as 
\[
EN=\frac{d}{2}E\left( \frac{1}{t\left(1-t\right)}-1 \right)\]
with 
\[
E\frac{1}{t\left(1-t\right)}=\frac{B\left(\nu-1,\mu-1\right)}{B\left(\nu,\mu\right)}.
\]
The second order moment is similarly computed as
\[
EN^{2}=\frac{d}{2}E\frac{\left(1-t\left(1-t\right)\right)\left(1+\frac{d}{2}\left(1-t\left(1-t\right)\right)\right)}{t^{2}\left(1-t\right)^{2}}=\frac{d}{2}\left(1+ \frac{d}{2} \right)E \frac{1}{t^{2}\left(1-t\right)^{2}}
-\frac{d\left( 1+d \right)}{2}E \frac{1}{t\left(1-t\right)} 
+\frac{d^{2}}{4} 
\]
so that the variance reads
\[
\sigma_{N}^{2}=\frac{d^{2}}{4}-\frac{d\left( 1+d \right)}{2}\frac{B\left(\nu-1,\mu-1\right)}{B\left(\nu,\mu\right)}
+\frac{d}{2}\left( 1+\frac{d}{2} \right)\frac{B\left(\nu-2,\mu-2\right)}{B\left(\nu,\mu\right)}-\frac{d^{2}}{4}\left(\frac{B\left(\nu-1,\mu-1\right)}{B\left(\nu,\mu\right)}-1\right)^{2}
\]
and the result follows after simplification.

\subsection{Proof of Theorem \ref{thm:thm3}}
It is clearly enough to prove that $\lim_{\nu\to\infty}\alpha_0^{(\nu,\gamma)}=1$ since
$0\le \alpha_k^{(\nu,\gamma)}\le 1-\alpha_0^{(\nu,\gamma)}$
for $k>0$. We have
\[
\alpha_0^{(\nu,\gamma)}=\frac{1}{\Gamma(\nu)}\int_0^\infty e^{-a}a^{\nu-1}
\left(1+\frac{\gamma^2}{a}\right)^{-d/2}\,da,
\]
hence
\[
1-\alpha_0^{(\nu,\gamma)}=\frac{1}{\Gamma(\nu)}\int_0^\infty e^{-a}a^{\nu-1}
\left(1-\left(1+\frac{\gamma^2}{a}\right)^{-d/2}\right)\,da=
\frac{d\gamma}{\Gamma(\nu)}\int_0^\infty e^{-a}a^{\nu-2}
\xi\left(1+\frac{\xi^2}{a}\right)^{-d/2-1}\,da,
\]
for some $0<\xi<\gamma$ by the mean value theorem. For $\nu>1$ we therefore get
\[
1-\alpha_0^{(\nu,\gamma)}=\frac{d\gamma^2}{\Gamma(\nu)}\int_0^\infty
e^{-a}a^{\nu-2}\,da=\frac{d\gamma^2}{\nu-1}
\]
which tends to 0 for $\nu\to\infty$.

\subsection{Proof of Theorem \ref{thm:thm4}}

The coefficients $\alpha_{k}^{(\nu,\tilde{\gamma})}$ read
\begin{eqnarray*}
\alpha_{k}^{(\nu,\tilde{\gamma})} & = & \frac{\Gamma(k+\frac{d}{2})}{k!\Gamma\left(\frac{d}{2}\right)\Gamma(\nu)}\left(\frac{\eta^{2}\left(\nu-1\right)}{\sigma^{2}}\right)^{k}\int_{0}^{+\infty}\exp(-a)a^{\nu+\frac{d}{2}-1}\left(a+\frac{\eta^{2}\left(\nu-1\right)}{\sigma^{2}}\right)^{-k-\frac{d}{2}}da\\
 & = & \frac{\Gamma(k+\frac{d}{2})}{k!\Gamma\left(\frac{d}{2}\right)\Gamma(\nu)}\left(\frac{\eta^{2}}{\sigma^{2}}\right)^{k}\left(\nu-1\right)^{\nu}\int_{0}^{+\infty}\exp(-b\left(\nu-1\right))b^{\nu+\frac{d}{2}-1}\left(b+\frac{\eta^{2}}{\sigma^{2}}\right)^{-k-\frac{d}{2}}db\end{eqnarray*}
using change of variable $a=b\left(\nu-1\right).$ Let us apply  Laplace's method \cite[p.278]{Widder}: if $f'\left( x_{0} \right)=0$ and $f''\left( x_{0} \right)<0$ then for large $M$,
\[
\int_{a}^{b} g(x) e^{Mf(x)} dx \sim g\left( x_{0} \right)e^{Mf\left( x_{0} \right)}\sqrt{\frac{2\pi}{M \vert f''\left( x_0 \right)\vert }}.
\]
Choosing $M=\nu -1,\, f\left( x \right) = -x + \log{x},\, g(x)=x^{\frac{d}{2}}\left(x+a^2  \right)^{-k-\frac{d}{2}}$ and $x_{0}=1$, 
an equivalent of the right-hand side integral reads
\[
e^{-\left(\nu-1\right)}\sqrt{\frac{2\pi}{\nu-1}}\left(1+\frac{\eta^{2}}{\sigma^{2}}\right)^{-k-\frac{d}{2}}.
\]
Thus an equivalent of the coefficients $\alpha_{k}^{\left(\nu,\tilde{\gamma}\right)}$
reads
\[
\frac{\Gamma(k+\frac{d}{2})}{k!\Gamma\left(\frac{d}{2}\right)}\frac{\left(\nu-1\right)^{\nu-\frac{1}{2}}e^{-\left(\nu-1\right)}\sqrt{2\pi}}{\Gamma\left(\nu\right)}\left(\frac{\eta^{2}}{\sigma^{2}}\right)^{k}\left(1+\frac{\eta^{2}}{\sigma^{2}}\right)^{-k-\frac{d}{2}}
\]
which, by Stirling's formula, is equivalent to
\[
\frac{\Gamma\left(k+\frac{d}{2}\right)}{k!\Gamma\left(\frac{d}{2}\right)}\left(\frac{\eta^{2}}{\sigma^{2}}\right)^{k}\left(1+\frac{\eta^{2}}{\sigma^{2}}\right)^{-k-\frac{d}{2}}.
\]
As a consequence, for large values of $\nu$, the density $f_{\mathbf{Z}}$
of $\mathbf{Z}$ reads
\begin{eqnarray*}
f_{\mathbf{Z}}\left(\mathbf{z}\right) & \sim & \sum_{k=0}^{+\infty}\frac{\Gamma\left( \frac{d}{2} \right)}{\left( \sigma\sqrt{2\pi} \right)^{d}\Gamma\left(k+\frac{d}{2}\right)}\left(\frac{\Vert\mathbf{z}\Vert^{2}}{2\sigma^{2}}\right)^{k}\exp\left(-\frac{\Vert\mathbf{z}\Vert^{2}}{2\sigma^{2}}\right)\frac{\Gamma\left(k+\frac{d}{2}\right)}{k!\Gamma\left(\frac{d}{2}\right)}\left(\frac{\eta^{2}}{\sigma^{2}}\right)^{k}\left(1+\frac{\eta^{2}}{\sigma^{2}}\right)^{-k-\frac{d}{2}}\\
 & = & \frac{1}{\left( \sigma\sqrt{2\pi} \right)^{d}}\frac{1}{\left( 1+\frac{\eta^{2}}{\sigma^{2}} \right)^{\frac{d}{2}}}\exp\left(-\frac{\Vert\mathbf{z}\Vert^{2}}{2\sigma^{2}}\right)\exp\left(\frac{\Vert\mathbf{z}\Vert^{2}}{2\sigma^{2}}\frac{\frac{\eta^{2}}{\sigma^{2}}}{1+\frac{\eta^{2}}{\sigma^{2}}}\right)\\
 & = & \frac{1}{\left( \tilde{\sigma}\sqrt{2\pi} \right)^{d}}\exp\left(-\frac{\Vert\mathbf{z}\Vert^{2}}{2\tilde{\sigma}^{2}}\right)
\end{eqnarray*}
with
\[
\tilde{\sigma}^{2}=\sigma^{2}+\eta^{2}.
\]

\subsection{Proof of Proposition \ref{pro:Student+Gaussian}}

Conditionnally to a given value of $K=k,$ the random variable $\sigma\Vert\mathbf{X}_{2k+d}\Vert \mathbf{U}_{d}$
follows the distribution $g_{k,\sigma}$; let us check that the random
variable $K$ follows a $NB\left(\frac{d}{2},\frac{a}{1+a}\right)$
compound distribution where $a$ is $\Gamma\left(\nu,\gamma^{2}\right)$:
conditionnally to $a,$ such a compound distribution reads\[
p_{k}\vert a=\binom{k+\frac{d}{2}-1}{k}\left(\frac{a}{1+a}\right)^{\frac{d}{2}}\left(\frac{1}{1+a}\right)^{k}\]
so that 
\begin{eqnarray*}
\alpha_{k}^{\left(\nu,\gamma\right)}=Ep_{k}\vert a & = & \binom{k+\frac{d}{2}-1}{k}\int_{0}^{+\infty}\frac{\gamma^{2\nu}}{\Gamma\left(\nu\right)}e^{-\gamma^{2}a}a^{\nu+\frac{d}{2}-1}\left(1+a\right)^{-k-\frac{d}{2}}da\\
 & = & \frac{\Gamma\left(k+\frac{d}{2}\right)\gamma^{2k}}{k!\Gamma\left( \frac{d}{2} \right)\Gamma\left(\nu\right)}\int_{0}^{+\infty}e^{-b}b^{\nu+\frac{d}{2}-1}\left(\gamma^{2}+b\right)^{-k-\frac{d}{2}}db
\end{eqnarray*}
by the change of variable $b=\gamma^{2}a;$ this expression coincides
with (\ref{eq:alphakmultivariate}).

\subsection{Proof of Proposition \ref{pro:Gaussian+Gaussian}}

Conditionnally to a value of $N=n,$ the random variable $\Vert\mathbf{T}_{2n+d}\Vert \mathbf{U}_{d}$
follows the distribution $\phi_{n,\nu+\mu};$ we check that the random
variable $N$ follows a $NB\left(\frac{d}{2},t\left(1-t\right)\right),$
where $t$ follows a Beta distribution with parameters $\mu$ and
$\nu$: conditionnally to $t,$ this compound distribution reads\[
p_{n}\vert t=\binom{n+\frac{d}{2}-1}{n}t^{\frac{d}{2}}\left(1-t\right)^{\frac{d}{2}}\left(1-t\left(1-t\right)\right)^{n}\]
so that 
\begin{eqnarray*}
c_{n}^{\left(\nu,\mu\right)} & = & Ep_{n}\vert t=\frac{1}{B\left(\nu,\mu\right)}\binom{n+\frac{d}{2}-1}{n}\int_{0}^{1}t^{\mu+\frac{d}{2}-1}\left(1-t\right)^{\nu+\frac{d}{2}-1}\left(1-t\left(1-t\right)\right)^{n}dt\\
 & = & \frac{1}{B\left(\nu,\mu\right)}\frac{\Gamma\left(n+\frac{d}{2}\right)}{n!\Gamma\left(\frac{d}{2}\right)}\int_{0}^{1}t^{\mu+\frac{d}{2}-1}\left(1-t\right)^{\nu+\frac{d}{2}-1}\left(1-t\left(1-t\right)\right)^{n}dt
\end{eqnarray*}
which coincides with (\ref{eq:cn}).

\section*{Acknowledgements}
The work of the first author has been supported by grant 272-07-0321 from the Danish Research Council for Nature and Universe. CV thanks CB for his invitation to the University of Copenhagen in February 2009.
The paper was finished during a stay of the first author as Guest Professor at University of Marne la Vall\'{e}e in May 2009.


\begin{thebibliography}{7}

\bibitem{Abramowitz}M. Abramowitz, I. Stegun, Handbook of Mathematical
Functions: with formulas, Graphs, and Mathematical Tables, Dover Publications,
1965.

\bibitem{Andrews2}G. Andrews, Euler's 'exemplum memorabile inductionis
fallacis' and $q-$Trinomial Coefficients, J. Amer. Math. Soc. 3,
653--669, 1990.

\bibitem{Andrews}G. E. Andrews, R. Askey, and R. Roy, Special Functions,
Cambridge University Press, Cambridge, UK, 1st ed., Feb. 1999.

\bibitem{Bening}V. E. Bening, V. Yu. Korolev, On an application of the
  Student distribution in the theory of probability and mathematical
  statistics, Theory Probab. Appl., Vol. 49, No. 3 (2005), 377--391.

\bibitem{Berg} C. Berg,  On powers of Stieltjes moment sequences, II,
J. Comput. Appl. Math., Vol. 199 (2007), 23--38.

\bibitem{BergVignat1}C. Berg, C. Vignat, Linearization Coefficients of Bessel Polynomials and Properties of Student t-Distributions, Constructive Approximation, Vol. 27 (2008), 15--32.

\bibitem{BergVignat2}C. Berg, C. Vignat, On some results of Cufaro Petroni about Student t-processes, J. Phys. A: Math. Theor. 41 (2008) 265004.

\bibitem{Cufaro}N. Cufaro Petroni, Mixtures in nonstable L\'{e}vy processes, J. Phys. A: Math. Theor. 40 (2007), 2227--2250.

\bibitem{Fang}K.-T. Fang, S. Kotz, K. W. Ng, Symmetric Multivariate
and Related Distributions, Monographs on Statistics and Applied Probability,
Chapman \& Hall/CRC, 1989.

\bibitem{Gradshteyn}I. S. Gradshteyn, I. M. Ryzhik,  Table of Integrals,
Series, and Products, Seventh Edition, Academic Press, 2007.

\bibitem{Leonenko}C. C. Heyde, N. N. Leonenko, Student Processes,
  Adv. Appl. Prob. 37 (2005), 342--365.

\bibitem{Keilson}J. Keilson, F.W. Steutel, Mixtures of Distributions,
Moment Inequalities and Measures of Exponentiality and Normality,
Annals of Probability, Volume 2, Issue 1, Feb. 1974, 112--130. 

\bibitem{Nason}G. P. Nason, On the sum of $t$ and Gaussian random
variables, Statistics and Probability Letters 76 (2006) 1280.

\bibitem{Ruben}H. Ruben,  On the distribution of the weighted
  difference of two independent Student
  variables. J. Roy. Statist. Soc. Ser. B, 22 (1960), 188--194.

\bibitem{Walker} G.A. Walker, J.G. Saw, The Distribution of Linear Combinations of $t$-Variables. Journal of the American Statistical Association, Vol. 73, Issue 364, Dec. 1978, 876--878. 

\bibitem{Widder}D.V. Widder, The Laplace Transform, Princeton University Press, 1946.

\bibitem{Witkovsky}V. Witkovsk\'{y}, Exact distribution of positive linear combinations of inverted chi-square random variables with odd degrees of freedom. Statistics and Probability Letters 56 (2002), 45--50. 
 
\end{thebibliography}
\end{document}